\documentclass[11pt]{amsart}
\usepackage{amsmath}
\usepackage{amsxtra}
\usepackage{amscd}
\usepackage{amsthm}
\usepackage{amsfonts}
\usepackage{amssymb}
\usepackage{psfig}
\usepackage{eucal}
\newcommand{\psdraw}[3]{\begin{array}{c} \hspace{-1mm}
\raisebox{-4pt}{\psfig{figure=#1.ps,width=#2,height=#3}}
\hspace{-1mm}\end{array}}
\newtheorem{thm} {Theorem}

\newtheorem{conj} {Conjecture}
\newtheorem{cor} {Corollary}
\newtheorem{prop} {Proposition}
\newtheorem{lem} {Lemma}
\newtheorem{defn} {Definition}

\title{A New Class of Wilf-Equivalent Permutations} 

\author{Zvezdelina Stankova-Frenkel and Julian West}
\address{Zvezdelina Stankova-Frenkel, Dept. of Mathematics and
Computer Science, Mills College, Oakland, CA, {\tt stankova@mills.edu}}
\address{Julian West, Dept. of Mathematics and Statistics, University
of Victoria, Canada, {\tt westj@mala.bc.ca}}

\begin{document}
\thispagestyle{empty}
\maketitle
\vspace*{-5mm}
\centerline{March 2001}

\begin{abstract}
For about 10 years, the classification of permutation patterns was
thought completed up to length 6. In this paper, we establish a new
class of Wilf-equivalent permutation patterns, namely,
$(n-1,n-2,n,\tau)\sim (n-2,n,n-1,\tau)$ for any $\tau\in S_{n-3}$.  In
particular, at level $n=6$, this result includes the only missing
equivalence $(546213)\sim(465213)$, and for $n=7$ it completes the
classification of permutation patterns by settling all remaining cases
in $S_7$.
\end{abstract}

\tableofcontents
\section{Introduction}

A permutation $\tau$ of length $k$ is written as
$(a_1,a_2,\ldots,a_k)$ where $\tau(i)=a_i,\,1\leq i \leq k$. For $k<10$
we suppress the commas without causing confusion. As usual, $S_n$
denotes the symmetric group on $[n]=\{1,2, ...,n\}$.

\begin{defn}
{\rm Let $\tau$ and $\pi$ be two permutations of lengths $k$ and $n$,
respectively. We say that $\pi$ is $\tau$-{\it avoiding} if there is no
subsequence $i_{\tau(1)}, i_{\tau(2)}, ..., i_{\tau(k)}$ of
$[n]$ such that $\pi(i_1)<\pi(i_2)<\ldots<\pi(i_k)$.
If there is such a subsequence, we say that the subsequence
$\pi(i_{\tau(1)}),$ $\pi(i_{\tau(2)})$,..., $\pi(i_{\tau(k)})$ is {\it of
type} $\tau$.}
\end {defn}

\smallskip
\noindent For example, the permutation $\omega=(52687431)$ avoids
$(2413)$ but does not avoid $(3142)$ because of its subsequence
$(5283)$.  An equivalent, but perhaps more insightful, definition is
the following reformulation in terms of matrices.

\smallskip
\begin{defn}
{\rm Let $\tau \in S_n$. The {\it permutation matrix} $M(\tau)$ is the
$n\times n$ matrix having a $1$ in position $(i,\tau(i))$ for
$1\leq i \leq n$, and having $0$ elsewhere.\footnote{To keep the
resemblance with the ``shape'' of $\tau$, we coordinatize $M(\tau)$
from the bottom left corner.} Given two permutation matrices $M$ and
$N$, we say that $M$ {\it avoids} $N$ if no submatrix of $M$ is identical
to $N$.}
\end{defn}

Note that a permutation matrix $M$ of size $n$ is simply a transversal
of an $n\times n$ matrix, i.e. an arrangement of 1's for which there is
exactly one 1 in every row and in every column of $M$. It is clear
that a permutation $\pi \in S_n$ contains a subsequence $\tau \in S_k$
if and only if $M(\pi)$ contains $M(\tau)$ as a submatrix.

\medskip
Let $S_n(\tau)$ denote the set of $\tau$-avoiding permutations in
$S_n$. 

\begin{defn} {\rm Two permutations $\tau$ and $\sigma$ are called {\it
Wilf-equivalent} if they are equally restrictive:
$|S_n(\tau)|=|S_n(\sigma)|$ for all $n \in N$. We denote this by $\tau
\sim \sigma$. If $|S_k(\tau)|=|S_k(\sigma)|$ for $k\leq n$, then we say
that $\tau$ and $\sigma$ are {\it equinumerant} up to level $n$.}
\end{defn} 

The basic problem in the theory of forbidden subsequences is to
classify all permutations up to Wilf-equivalence. Obviously, if two
permutations are Wilf-equivalent, then they must be of the same
length. Further, many Wilf-equivalences can be deduced by symmetry
arguments within the same $S_k$.  For instance, if $M(\pi)$ contains
$M(\tau)$ as a submatrix, then the transpose matrix
$M(\pi)^t=M(\pi^t)$ contains $M(\tau)^t=M(\tau^t)$. The same is true
when simultaneously reflecting both matrices $M(\pi)$ and $M(\tau)$ in
either a horizontal or a vertical axis of symmetry. The three
operations defined above generate the dihedral group $D_4$ acting on
the set of permutation matrices in the obvious way. The orbits of
$D_4$ in $S_k$ are called {\it symmetry classes}. It is clear that if
$\tau$ and $\sigma$ belong to the same symmetry class in $S_k$, then
$\tau \sim \sigma$. However, Wilf-classes are in general, but
apparently rarely, larger than single symmetry classes. This makes the
classification of permutations up to Wilf-equivalence a subtle and
difficult process.

The first major result in the theory of forbidden subsequences states
that $(123)\sim (132)$, and hence $S_3$ is one Wilf-class, which
combines the two symmetry classes of $(123)$ and $(132)$.  As the
behest of Wilf, bijections between $S_n(123)$ and $S_n(132)$ were
given by Simion-Schmidt \cite{SS}, Rotem \cite{Ro}, Richards
\cite{Ri}, and West \cite{We1}. They all prove $|S_n(123)|=c_n$, where
$c_n$ is the $n$th Catalan number.  Permutations with forbidden
subsequences arise naturally in computer science in connection with
sorting problems and strings with forbidden subwords.  For example, in
\cite{Kn1}-\cite{Kn2} Knuth shows that $S_n(231)$ is the set of
stack-sortable permutations (see also \cite{Lov}), so that
$|S_n(231)|$ is the number of binary strings of length $2n$, in which
0 stands for a ``move into a stack'' and 1 symbolizes a ``move out
from the stack''.

Numerous problems involving forbidden subsequences have also appeared
in algebraic combinatorics.  In the late 1980s, it was discovered that
the property of avoiding 2143 exactly characterizes the {\it
vexillary} permutations, i.e. those whose Stanley symmetric function
is a Schur function. (See \cite{Macdonald} for a good exposition.)
Lakshimibai and Sandhya \cite{Lak} likewise show that $S_n(3412,4231)$
is the set of permutations indexing an interesting subclass of
Schubert varieties.  And Billey and Warrington \cite{Billey} have very
recently defined a class of permutations under 5 restrictions which
are related to the Kazhdan-Lusztig polynomials. This all naturally
leads to the study and classification of Wilf-classes of permutations
of length 4 or more.

\smallskip
The classification of $S_4$ turns out to be much more complicated than
that of $S_3$. It is completed in a series of papers by Stankova and
West. They utilize the concept of a {\it generating tree} ${\mathfrak
T}(\tau)$ of $\tau\in S_k$: the nodes on the $k$th level of
${\mathfrak T}(\tau)$ are the permutations in $S_n(\tau)$, and the
descendants of $\pi\in S_n(\tau)$ are obtained from $\pi$ by inserting
$n+1$ in appropriate places in $\pi$. Clearly, the tree isomorphism
${\mathfrak T}(\tau)\simeq {\mathfrak T}(\sigma)$ implies $\tau\sim
\sigma$, but the converse is far from true. In \cite{We1}, West shows
${\mathfrak T}(1234)\simeq {\mathfrak T}(1243)\simeq {\mathfrak
T}(2143)$. In \cite{St1}, Stankova constructs a specific isomorphism
${\mathfrak T}(4132)\cong {\mathfrak T}(3142)$. In \cite{St2}, she
completes the classification of $S_4$ by proving $(1234)\sim(4123)$;
there she uses a different approach which yields the somewhat
surprising result that, while ${\mathfrak T}(1234)\not\cong{\mathfrak
T}(4123)$, on every level of the two trees the number of nodes with a
given number of descendants is the same for both trees.  Thus, the
seven symmetry classes of $S_4$ are grouped in three Wilf-classes,
with representatives (4132), (1234) and (1324)
(cf. Fig.~\ref{fig-S4}.)

\smallskip
In \cite{BW}, Babson-West show $(n-1,n,\tau)\sim(n,n-1,\tau)$ for any
$\tau\in S_{n-2}$, and $(n-2,n-1,n,\tau)\sim(n,n-1,n-2,\tau)$ for any
$\tau\in S_{n-3}$, thus completing the classification up to level
5. The key idea is the concept of a stronger Wilf-equivalence
relation.

\begin{defn}{\rm  A {\it transversal} $T$ of a Young
diagram $Y$ is an arrangement of 1's and 0's such that every row and
every column of $Y$ has exactly one 1 in it. A subset of 1's in $T$ is
said to form a {\it submatrix} of $Y$ if all columns and rows of $Y$
passing through these 1's intersect {\it inside} $Y$. For a
permutation $\tau\in S_k$, we say that $T$ {\it contains} pattern
$\tau$ if some $k$ 1's of $T$ form a submatrix of $Y$ identical to
$M(\tau)$ (cf. Fig.~\ref{fig-splitting}.)}
\end{defn}

Given several 1's in a transversal $T$, the condition for them to form
a submatrix of $T$ is the same as the requirement that the column of
the rightmost 1 and the row of the lowest 1 must intersect {\it
inside} $Y$. This condition is necessary for the new definition to be
a useful generalization of the classical definition of a forbidden
subsequence, as we shall see below.  In particular, when $Y$ is a
square diagram, the two definitions coincide. Let us denote by
$S_Y(\tau)$ the set of all transversals of $Y$ which avoid $\tau$.

\begin{defn}{\rm Two permutations $\tau$ and $\sigma$
are called {\it shape-Wilf-equivalent} (SWE) if
$|S_{Y}(\tau)|=|S_{Y}(\sigma)|$ for all Young diagrams $Y$.  We denote
this by $\tau\stackrel{s}{\sim}\sigma$.}
\end{defn}

Clearly, $\tau\stackrel{s}{\sim}\sigma$ implies $\tau{\sim}\sigma$,
but not conversely.  We will write $Y(a_1,a_2,...,a_n)$ for the Young
diagram $Y$ whose $i$-th row has $a_i$ cells, $1\leq i\leq
n$.  In order for a Young diagram $Y$ to have any transversals at all,
$Y$ must have the same number of rows and columns and $Y$ must contain
the {\it staircase} diagram $St=Y(n,n-1,...,2,1)$, where $n$ is the
number of cells in the top (largest) row of $Y$. Thus, from now on,
when we talk about a Young diagram $Y$ of size $n$, we will assume
that $Y$ has $n$ rows and $n$ columns and contains $St$ of size $n$.

SWE is a very strong relation on two permutations and it is certainly too
restrictive on its own to be useful in the general classification of
permutations. However, combined with the proposition below
(see \cite{BW}), it allows for more Wilf-equivalences to be
established.

\begin{prop}  Let $A\stackrel{s}{\sim}B$ for some permutation
matrices $A$ and $B$. Then for any permutation matrix $C$:
\[\left( \begin{array}{c|c}
          A & 0\\ \hline
          0 & C
      \end{array} \right) \stackrel{s}{\sim}
\left(\begin{array}{c|c}
          B & 0\\ \hline
          0 & C
      \end{array} \right)\cdot\]
\label{prop-BW}
\end{prop}

Let $I_k$ be the $k\times k$ identity matrix, and let $J_k$ be its
reflection across a vertical axis of symmetry.  According to
Backelin-West-Xin in \cite{BWX}, $I_k\stackrel{s}{\sim} J_k$ for any
$k$, and hence $(n,n-1,...,m,\tau)\sim (m,...,n-1,n,\tau)$ for any
$\tau\in S_{n-m}$. This SWE generalizes the results in \cite{We2} and
\cite{BW}, but it is not sufficient to complete the classification of
$S_7$, nor of $S_6$.

\begin{figure}[h]
{\small \begin{tabular}{|c||c|c|c|c|c|c|c|}
\hline\hline
$\tau\in S_6$&$S_8(\tau)$&$S_9(\tau)$
&$S_{10}(\tau)$&$S_{11}(\tau)$&$S_{12}(\tau)$&$S_{13}(\tau)$\\\hline
\hline
(546213) & 39428 & 344772 & 3289163 & 33765743 & 368833207 & 4247979687 \\
\hline
(465213) & 39428 & 344772 & 3289163 & 33765743 & 368833207 & 4247979687
\\\hline
\end{tabular}}
\caption{Missing pair of Wilf-equivalence in $S_6$}
\end{figure}

In 2001, Stankova-Frenkel noticed a missing case of a plausible
Wilf-equivalence in $S_6$: (546213) and (465213) were equinumerant up
to level $11$, but no reference was found regarding why these
permutations were thought to be in different Wilf
classes. Stankova-Frenkel further found an infinite class of
Wilf-equivalences
\begin{equation}
(n-1,n-2,n,\tau)\sim (n-2,n,n-1,\tau),
\label{zvez} 
\end{equation}
At her request, West confirmed (\ref{zvez}) by computer checks for
$n=6,7$ up to level 13.\footnote{In all tables, we skip the column
corresponding to $|S_{n+1}(\tau)|$ if $\tau\in S_n$. It is easy to see
that all permutations in $S_n$ are equinumerant on level $n+1$ (for
example, cf. \cite{NW}).}

\smallskip The purpose of this paper is to explain the proof of the
new Wilf-equivalences in (\ref{zvez}). The idea is to show
$(213)\stackrel{s}{\sim} (132)$ and apply then
Proposition~\ref{prop-BW}. Even though $M(213)$ and $M(132)$ are
transposes of each other, their SWE relationship is far from trivial,
as the present paper will reveal.  It is surprising that, since
the launch of the theory of forbidden subsequences in the early 1980s,
and the introduction of SWE in the early 1990s, such a basic
relationship is discovered only now.

\begin{thm}[Main result of the paper]
The permutations $(213)$ and $(132)$ are
shape-Wilf-equivalent. Consequently, for any $\tau\in S_{n-3}$, the
permutations $(n-1,n-2,n,\tau)$ and $(n-2,n,n-1,\tau)$ are
Wilf-equivalent.
\label{thm-Wilf}
\end{thm}
\noindent Theorem 1 finally accounts for the last missing case in
$S_6$ and the remaining cases in $S_7$, thus completing the
classification of forbidden subsequences up to length $n=7$.

\begin{figure}[h]
{\small \begin{tabular}{|c||c|c|c|c|c|}
\hline\hline
$\tau\sim \sigma$ 
in $S_7$&$S_9(\tau)$&$S_{10}(\tau)$&$S_{11}(\tau)$&$S_{12}(\tau)$
& $S_{13}(\tau)$\\\hline
\hline
(6571342)$\sim$(5761342)
& 361300 & 3587768 & 38951398 & 457416920 & 5756026177 \\ \hline
(6571423)$\sim$(5761423)
& 361300 & 3587768 & 38951411 & 457418106 & 5756088993 \\ \hline
(6572413)$\sim$(5762413)
& 361300 & 3587768 & 38951430 & 457419793 & 5756176230 \\ \hline
(6572431)$\sim$(5762431)
& 361300 & 3587768 & 38951467 & 457423216 & 5756360170 \\ \hline
(6574132)$\sim$(5764132)
& 361300 & 3587780 & 38952330 & 457459680 & 5757549454 \\ \hline
(6571432)$\sim$(5761432)
& 361301 & 3587834 & 38953996 & 457496956 & 5758168203 \\ \hline
(6571243)$\sim$(5761243)
& 361301 & 3587834 & 38954024 & 457499462 & 5758298471 \\ \hline
\end{tabular}}
\caption{Final Wilf-equivalences in $S_7$}
\label{new-Wilf-S7}
\end{figure}

\smallskip
In summary, modulo symmetry classes, as of now there are essentially
two known infinite families of Wilf-equivalences, resulting from
\cite{BWX} and the present paper:
\begin{equation}
\left(\begin{array}{c|c}
               I_k& 0\\ \hline
               0  & C
        \end{array}\right) \stackrel{s}{\sim}
\left(\begin{array}{c|c}
               J_k& 0\\ \hline
               0  & C
        \end{array}\right)\,\,\text{and}\,\,
\left(\begin{array}{c|c}
               M(213)& 0\\\hline
               0  & C
        \end{array}\right)\stackrel{s}{\sim}
\left(\begin{array}{c|c}
               M(132)& 0\\\hline
               0  & C
        \end{array}\right)\cdot
\label{infinite}
\end{equation}
Further, there is only one known ``sporadic'' case of Wilf-equivalence,
from \cite{St1}:
\begin{equation}
\left(\begin{array}{cccc}
               1 & 0 & 0 & 0\\
               0 & 0 & 1 & 0\\
               0 & 0 & 0 & 1\\
               0 & 1 & 0 & 0
        \end{array}\right) {\sim}
\left(\begin{array}{cccc}
               0 & 0 & 1 & 0\\
               1 & 0 & 0 & 0\\
               0 & 0 & 0 & 1\\
               0 & 1 & 0 & 0
        \end{array}\right)\cdot
\label{sporadic}
\end{equation}
Of all symmetry classes in $S_4$, the classes of the above (4132) and
(3142) in (\ref{sporadic}) are the most unexpected to fall into the
same Wilf-class: (3142) has the smallest symmetry class as it
corresponds geometrically to the quadrilateral with most symmetries -
the square, while $(4132)$ has the largest symmetry class as it
corresponds to the quadrilateral with least number of symmetries - a
quadrilateral with 4 different angles. And yet, not only $(4132)\sim
(3142)$, but also their trees are isomorphic. In a similar vein, the
permutations in (\ref{infinite}) are more than just Wilf-equivalent -
they are SWE. This is an interesting phenomenon - so far, every known
Wilf-equivalence can be explained by a stronger relationship: either
symmetry, tree isomorphism, SWE, or a combination of these.

\smallskip
For further discussion, we refer the reader to
Sections~\ref{section-results}-\ref{section-conjectures}.

\section{The Row-Decomposition Formula}

Let $Y$ be a Young diagram with $n$ rows and $n$ columns.  We denote
by $(k,l)$ the intersection cell of row $k$ and column $l$, counted
from the top left corner of $Y$. A cell in the bottom row of $Y$ is
called a {\it bottom} cell of $Y$. Let $m$ be the number of bottom
cells in $Y$.

\begin{defn}{\rm For a subset $X$ of cells in
$Y$, define the {\it reduction $Y\!\!\big/_{\displaystyle{\!\!X}}$ of
$Y$ along $X$} to be a new Young diagram obtained from $Y$ by deleting
all rows and columns of $Y$ which intersect $X$.}
\end{defn}

\smallskip For example, if the {\it cross} of a cell $C$ in $Y$ is the
union of the row and column containing $C$, then the reduction
$Y\!\!\big/_{\displaystyle{\!C}}$ is the diagram obtained from $Y$ by
deleting the cross of $C$. The reduction of $Y$ along an arbitrary
bottom cell of $Y$ is denoted by $Y^r$, and it clearly does not depend
on the choice of the bottom cell. We will use this fact frequently
when reducing along cell $(1,n)$ (and $(n,1)$) in the proof of the
commitativity argument in Lemma~\ref{lem-commute} in
Section~\ref{section-commute}.

\begin{defn}{\rm To any bottom cell $C$ in
$Y$, we associate a {\it cross-product} $Y_C$ of Young diagrams in the
following way. Mark by $P$ the top right corner of $C$ (this is a grid
point on $Y$), and consider the reduction
$Y\!\!\big/_{\displaystyle{\!C}}$, which still contains point
$P$. Starting from $P$ in $Y\!\!\big/_{\displaystyle{\!C}}$, draw a
$45^{\circ}$ ray in north-east direction until the ray intersects for
the first time the border of $Y\!\!\big/_{\displaystyle{\!C}}$, and
use the resulting segment as the diagonal of a smaller subdiagram
$A_C$ of $Y\!\!\big/_{\displaystyle{\!C}}$. Delete the rows and
columns of $A_C$ in $Y\!\!\big/_{\displaystyle{\!C}}$, leaving a
subdiagram $B_C=Y\!\!\big/_{\displaystyle{\!\{C,A_C\}}}$. Thus, any
bottom cell $C$ in the original diagram $Y$ determines a pair
$(A_C,B_C)$ of smaller Young diagrams, which we call the {\it
cross-product} of $A_C$ and $B_C$ and denote by $Y_C:=A_C\times B_C$.}
\end{defn}

\begin{figure}[h]
$$\psdraw{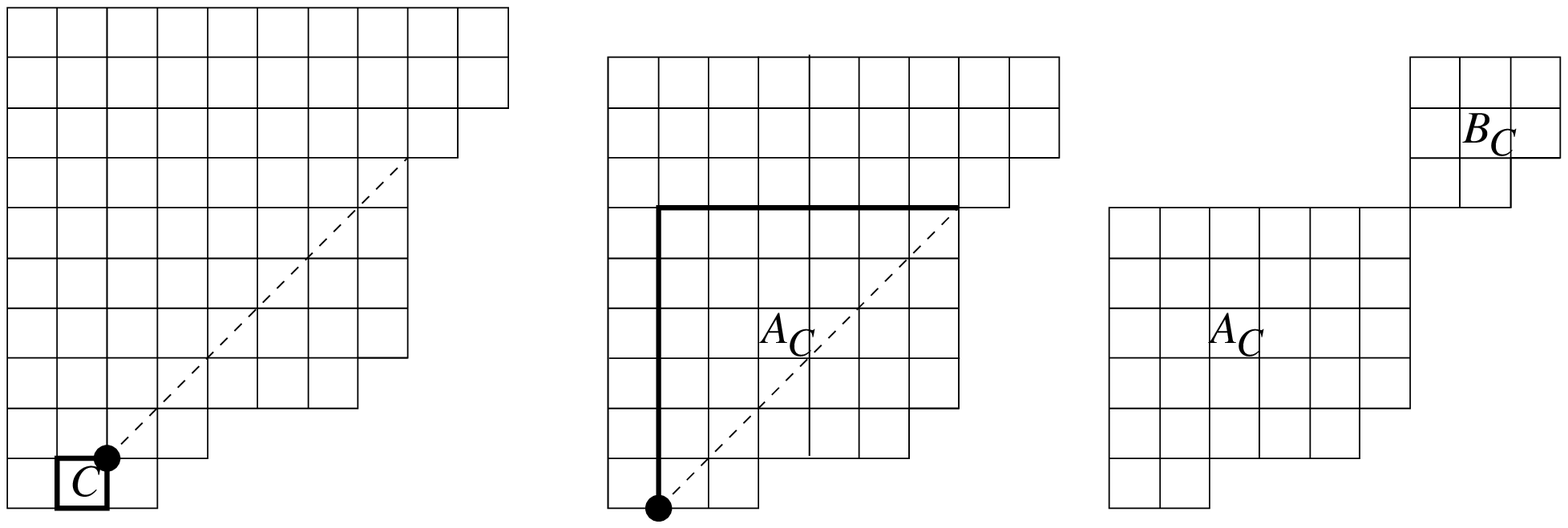}{3.2in}{1.1in}$$
\caption{$Y\rightarrow Y\!\!\big/_{\displaystyle{\!C}}\rightarrow Y_C=A_C\times
B_C$}
\label{cross-fig}
\end{figure}

\smallskip
If one of the subdiagrams $A_C$ or $B_C$ is empty, we define $Y_C$ to
equal the other subdiagram. This case occurs exactly when $C$ is the
first or the last bottom cell of $Y$:
\[Y_{(n,1)}
=Y^r\times \emptyset =Y^r=\emptyset \times Y^r=Y_{(n,m)}.\]
\noindent{\bf Example 1.} Let $Y=Y(10,10,9,8,8,8,8,7,4,3)$. Let
$C=(10,2)$. Then $Y_C=A_C\times B_C=Y(6,6,6,6,5,2)\times Y(3,3,2)$
(cf. Fig.~\ref{cross-fig}.)

\begin{defn}{\rm Let $Y$ be a Young diagram of size
$n$. The {\it row decomposition} of $Y$ is the formal sum ${\mathcal
R}(Y)$ of cross-products of smaller Young diagrams: \[{\mathcal
R}(Y):=\sum_C Y_{C}=\sum_C A_C\times B_C,\] where the sum is taken
over all bottom cells $C$ of $Y$.}
\end{defn}

\smallskip
\noindent As noted above, the first and the last summands of
${\mathcal R}(Y)$ are identical to $Y^r$.

\begin{defn}{\rm 
A transversal of a diagram $Y$ which avoids (resp.  contains) a
$(213)$--pattern is called a {\it good} (resp. {\it bad}) transversal
of $Y$. Denote by ${\mathcal T}(Y)={\mathcal T}(a_1,a_2,...,a_n)$ the
number of good transversals of $Y=Y(a_1,a_2,...,a_n)$.}
\end{defn}
We use the convention ${\mathcal T}(\emptyset)=1$. For some $1$
already placed in a cell of $Y$, we say that it {\it imposes a
$(213)$-condition on $Y$} if it plays the role of a ``1'' in a
$(213)$-pattern contained in some bad transversal of $Y$; the
(213)-condition is the actual condition on the rest of $Y$ in order to
avoid a (213)-pattern containing this 1.

\begin{defn}{\rm 
Two diagrams $Y$ and $X$ are said to be {\it numerically equivalent}
if ${\mathcal T}(Y)={\mathcal T}(X)$.  We denote this by $Y\equiv X$.}
\end{defn}

\smallskip Clearly, ${\mathcal T}(A_C\times B_C)={\mathcal
T}(A_C)\cdot {\mathcal T}(B_C)$. Moreover, to obtain the number
${\mathcal T}(Y)$, we can apply the function $\mathcal T$ to all terms
in the formal sum ${\mathcal R}(Y)$:

\begin{thm}[Row-Decomposition] Let $Y$ be a Young diagram of size $n$. Then 
\begin{equation}
{\mathcal T}(Y)=\sum_C {\mathcal T}(A_C)\cdot {\mathcal T}(B_C)
\label{RD-formula}
\end{equation}
where the sum is taken over all bottom cells $C$ in $Y$.
\label{thm-RD}
\end{thm}

\noindent
{\sc Proof:} For a bottom cell $C=(n,i)$, let $Y_i$ denote the diagram
$Y$ with the additional data of $1$ in cell $C$. We claim that the
good transversals of $Y_i$ are in 1-1 correspondence with the good
transversals of the pair of diagrams $(A_C,B_C)$, and hence
\begin{equation}
{\mathcal T}(Y_i)={\mathcal T}(A_C)\cdot {\mathcal T}(B_C).
\label{claim}
\end{equation} 
This, and the fact that any transversal of $Y$ must contain exactly
one 1 in the bottom row, immediately establishes
Theorem~\ref{thm-RD}. Hence, it suffices to prove (\ref{claim}) for
$Y_i$.

\smallskip
When $i=1$ or $i=m$, the claim (\ref{claim}) is trivial. Indeed, in
these cases, the 1 in the bottom row of $Y$ is either in the first or
last bottom cell, and hence it doesn't impose any (213)-conditions on
$Y$. The question reduces to finding all good transversals of
$Y\!\!\big/_{\displaystyle{\!C}}=Y^r=Y_C$. Therefore, $Y_1\equiv Y_m
\equiv Y^r$.

\smallskip
Assume now that $1<i<m$.  Fix a good transversal $T$ of $Y_i$. Denote
by $1_j$ the 1 in column $j$. By $1_j>1_k$ we mean that $1_j$ is in a
row above the row of $1_k$. Similarly, for two disjoint sets $A$ and
$B$ of 1's, by $A>B$ we mean that all 1's in A are above all 1's in B.
Let $B_L=\{1_1,1_2,...,1_{i-1}\}$ denote the set of all 1's in $T$
appearing in columns to the left of cell $C$. Similarly, let
$A=\{1_{i+1},1_{i+2},...,1_k\}$ be the set of all 1's appearing in the
columns of $Y$ intersecting $A_C$, and let
$B_R=\{1_{k+1},1_{k+2},...,1_n\}$ be the set of all 1's appearing in
the remaining columns, i.e. all columns to the right of $A_C$.  Notice
that no 1 in $B_R$ can appear in a row intersecting $A_C$: the rows of
the 1's in $B_R$ are above all rows of $A_C$ as enforced by the
construction of $A_C$ via the $45^{\circ}$ segment.

\smallskip The key idea of the proof is contained in the following
lemma:

\begin{lem} In any good transversal $T$ of $Y_i$, we have $B_L>A$. 
\label{lem-split}
\end{lem}

\noindent {\sc Proof} (of Lemma~\ref{lem-split}): Given a row $j$ of
$A$, let the {\it level} $L_j$ be the subset of $A$ consisting of all
1's whose orthogonal projections onto row $j$ are inside
$A_C$. Clearly, every $1\in A$ belongs to at least one level $L_j$,
and $L_{n-1}\subseteq L_{n-2}\subseteq \cdots \subseteq L_{n-k}=A$,
where $k$ is the size of $A_C$. Recall that $i<m$, so that
$L_{n-1}\cap A\not= \emptyset$, and that the bottom row of $Y$ is
filled by $1_i$. This imposes a (213)-condition on $Y_i$:
$B_L>L_{n-1}$.  Finally, from the construction of $A_C$, $|L_{n-j}|>j$
for $j=1,2,..,k-1$, and $|L_{n-k}|=k$.

\smallskip
We will prove simultaneously the following two statements for all
$1\leq j\leq k$:

\begin{itemize}
\item[(i)] $B_L>L_{n-j}$;

\item[(ii)] All rows $n-1$, $n-2$, ..., $n-j$ of $Y_i$ are filled in
with 1's on level $L_{n-j}$.
\end{itemize}

\smallskip 
For $j=1$, (i) was shown above. But then row $n-1$ in $Y_i$ must be
filled with an element of $L_{n-1}$, so (ii) is also true. Assume (i)
and (ii) for some $j<k$. Then (ii) for $j$, together with
$|L_{n-j}|>j$, implies that at least one element $1_s$ of $L_{n-j}$ is
in row $n-j-1$ or above. Now (i) for $j$ implies in particular
$B_L>1_s$, so that a new (213)-condition is imposed:
$B_L>(L_{n-j-1}\backslash L_{n-j})$. Combining, $B_L>L_{n-j-1}$: this
is (i). By definition of $L_{n-j-1}$, the only $1_h$ that can possibly
fill in row $n-j-1$ in $Y_i$ must belong either to $L_{n-j-1}$, or to
$B_L$. Because $|L_{n-j-1}|\geq j+1$ and $B_L>L_{n-j-1}$, we conclude
that $1_h\in L_{n-j-1}$. This shows (ii) for $j+1$ and completes the
inductive proof of the above statement. Lemma~\ref{lem-split} follows
automatically from (i) for $j=k$.\qed

\smallskip
\noindent{\sc End of Proof of Theorem~\ref{thm-RD}:} Combining Lemma
\ref{lem-split} with a previous observation, we see that no 1's from
$B_L$ or from $B_R$ can fill the rows intersecting $A_C$.  In other
words, all rows of $A_C$ must be filled exactly with the 1's from set
$A$: the number of necessary 1's to make a transversal of $A_C$
matches $|A|=k$ because, by construction, $A_C$ has as many rows as
columns. This in turn forces all 1's in $B_L$ and $B_R$ to make a good
transversal of the subdiagram $B_C$. It remains to show that there are
no further (213)-conditions imposed by a triple of 1's coming from
$A_C$ and $B_C$.

\smallskip
The only way for $A_C$ and $B_C$ to engage together in a (213)-pattern 
is to have the ``2'' in $B_L$, the ``1'' in $A$, and the ``3'' in $B_R$;
or, to have the ``2'' and the ``1'' in $A$, and the ``3'' in $B_R$. Even
though such configurations of three 1's are possible, their ``full''
matrices will not be contained entirely in $Y$ because of the relative
positioning of $A_C$ and $B_R$.

\smallskip
Putting everything together, the good transversals of $Y_i$ are in
1--1 correspondence with pairs of good transversals of $A_C$ and
$B_C$, i.e. ${\mathcal T}(Y_i)={\mathcal T}(A_C)\cdot {\mathcal T}(B_C)$ for all bottom cells $C$ of
$Y$. This completes the proof of the Row-Decomposition formula. \qed

\medskip
\noindent
{\bf Example 2.} To illustrate the above proof, consider
$Y(10,10,9,8,8,8,8,7,4,3)$. Let $Y_2$ denote the diagram $Y$ with the
additional data that the 1 in the bottom row is in cell $C=(10,2)$.
We have to show that ${\mathcal T}(Y_2)={\mathcal T}(6,6,6,6,5,2)\cdot
{\mathcal T}(3,3,2)$ (cf. Fig.~\ref{cross-fig}.)

\begin{figure}[h]
$$\psdraw{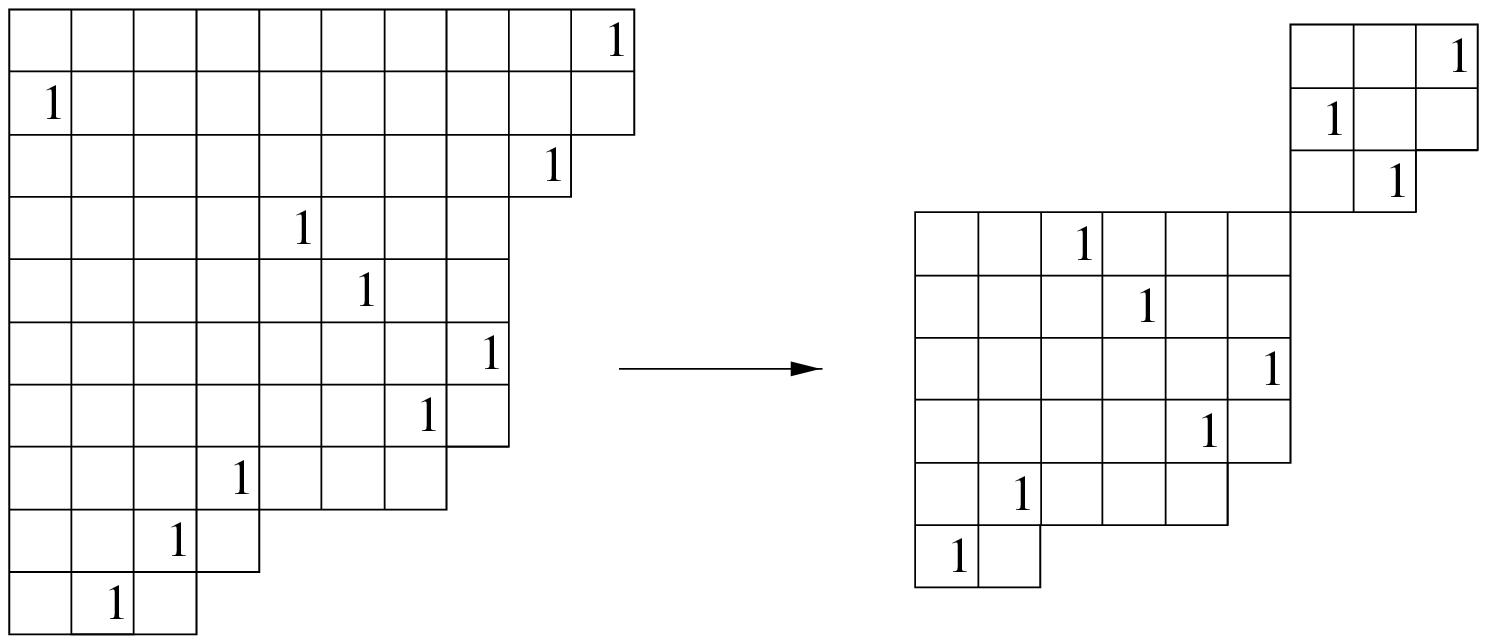}{2.5in}{1.1in}$$
\caption{Splitting of a transversal $T\in S_{Y(10,10,9,8,8,8,8,7,4,3)}
(213)$}
\label{fig-splitting}
\end{figure}

\smallskip
The initial condition that $1_2$ is in the bottom row forces the
(213)-condition $1_1> 1_3$. Since $1_3$ is above row 10, the only 1's
which can fill row 9 are $1_3$ and $1_4$.  If $1_3$ is in row 9, then
$1_1>1_4$ in order to avoid (213); if $1_4$ in in row 9, then
$1_1>1_3>1_4$. In any case, $1_1>1_3,1_4$. Without loss of generality,
assume that $1_3$ is in row 9, so that $1_4$ is in row 8 or
above. From $1_1>1_4$ and avoiding (213), we conclude that
$1_1>1_4,1_5,1_6,1_7$. One of the latter four 1's must fill in row
$8$. Without loss of generality, assume that $1_4$ is in row 8; hence
$1_5,1_6,1_7$ are in rows $7$ or above. But then, to avoid (213), we
are forced to conclude that $1_1>1_8$, i.e. $1_1$ is above all of
$1_3,...,1_8$. We need six 1's to fill in the six rows $9,8,...,4$. It
immediately follows that $1_3,...,1_8$ must have filled all 6 rows and
columns of subdiagram $A_C$, leaving all remaining 1's (except for
$1_2$) to form a good transversal of subdiagram $B_C$
(cf. Fig.~\ref{fig-splitting}.)

\smallskip
We note that $B_C$ consists of two disjoint parts: $B_L(1,1,1)$ and
$B_R(2,2,1)$. The argument that $B_C$ and $A_C$ cannot engage
together in a pattern (213) is identical to the corresponding
part of the proof of Theorem~\ref{thm-RD}. \qed

\medskip
Given a diagram $Y$, let $C_b$ be its rightmost bottom cell, which we
call the {\it bottom corner} of $Y$, and let $C_{b-1}$ be the bottom
cell to the left of $C_b$. (In our previous notation, $C_b=(n,m)$,
$C_{b-1}=(n,m-1)$.) Deleting $C_b$ from $Y$ results in a new diagram,
which we denote by $Y^{\prime}_r$ and call the {\it row-deletion of}
$Y$.  Similarly, we define the {\it right corner} $C_t$ of $Y$ as the
bottom cell in the rightmost column of $Y$, $C_{t-1}$ as the cell
directly above $C_t$, and the {\it column-deletion} $Y^{\prime}_c$ by
deleting $C_t$ from $Y$.  Note that
$(Y^{\prime}_r)^t=(Y^t)^{\prime}_c$, where $X^t$ denotes as usual the
transpose of diagram $X$ along its main (north-west to south-east)
diagonal.

\smallskip
In the row-decomposition of $Y$, we distinguish one special summand:
the last but one summand $Y_{C_{b-1}}=A_{C_{b-1}}\times B_{C_{b-1}}$,
which we denote by $Y^{\prime\prime}_r$. 

\begin{cor} For any diagram $Y$, $Y \equiv 
Y^{\prime}_r + Y^{\prime\prime}_r$.
\label{cor-RD}
\end{cor}

\noindent
{\sc Proof:} The row-decomposition of $Y$ includes one more summand
than the row-decomposition of $Y^{\prime}_r$, namely, $Y^{\prime\prime}_r$:
${\mathcal R}(Y)=Y^{\prime}_r+Y^{\prime\prime}_r$. 
Theorem~\ref{thm-RD} completes the proof.\qed

\medskip
From now on, we shall refer to Corollary~\ref{cor-RD} as
Row-Decomposition (RD).

\begin{figure}[h]
$$\psdraw{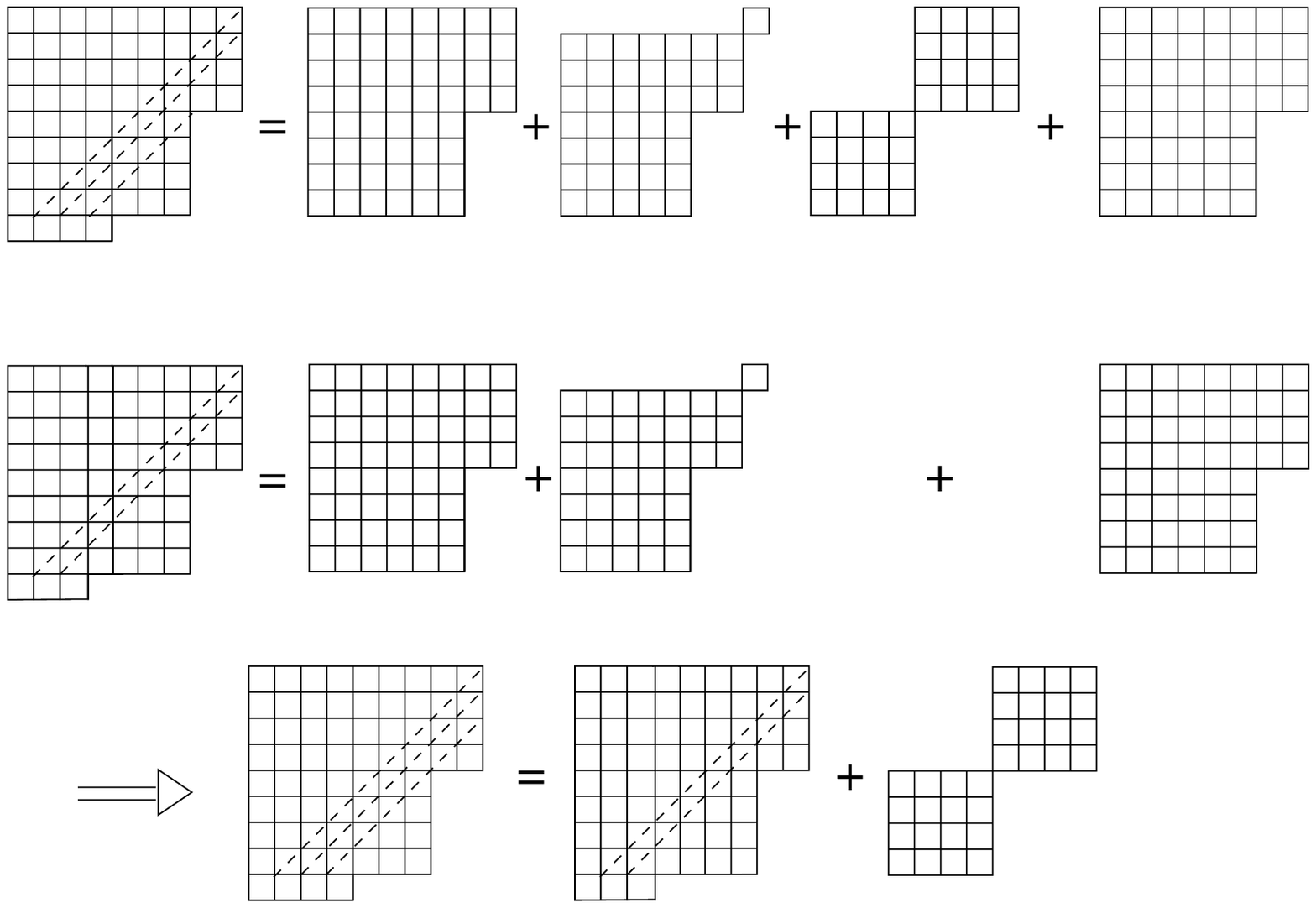}{4in}{2.9in}$$
\caption{${\mathcal T}(9,9,9,9,7,7,7,7,4)={\mathcal T}(9,9,9,9,7,7,7,7,3)+{\mathcal T}(4,4,4,4)^2$}
\label{fig-cor}
\end{figure}

\medskip
\noindent
{\bf Example 3.} Fig.~\ref{fig-cor} illustrates Corollary~\ref{cor-RD}
for ${\mathcal T}(9,9,9,9,7,7,7,7,4)$.

\section{Column Decomposition}

\begin{defn}{\rm
The {\it column decomposition} ${\mathcal C}(Y)$ of $Y$ is defined by:
\[{\mathcal C}(Y)=\sum_{C}((Y^t)_{C^t})^t\]
where $C$ runs over all cells in the rightmost column of $Y$, and
$C^t$ is the image of the cell $C$ after transposing $Y$.}
\end{defn}

\smallskip
\noindent Note that the column decomposition ${\mathcal C}(Y)$ can be obtained
directly from $Y$ without going through the transpose $Y^t$: for
a cell $C$ in the rightmost column of $Y$, mark the south-west corner
of $C$, draw the $45^{\circ}$ ray in south-west direction until it
intersects the border of $Y$, delete the cross of $C$ and define
analogously the product $Y_C=A_C\times B_C$.

\smallskip
As with row decomposition, we denote by $Y^{\prime}_c$ the diagram
resulting from $Y$ by deleting the right corner cell $C_t$, and by
$Y^{\prime\prime}_c$ the summand in ${\mathcal C}(Y)$ corresponding to the cell
$C_{t-1}$ right above $C_t$. By definition, it is clear that 
\[{\mathcal C}(Y)=Y^{\prime}_c+Y^{\prime\prime}_c.\]
It is not obvious, however, why the same formula should be true {\it
after} applying the function $\mathcal T$ to all terms: why is
$Y\equiv Y^{\prime}_c+Y^{\prime\prime}_c$?

\begin{lem} The following statements are equivalent:

\begin{itemize}
\item [(i)] $(213)\stackrel{S}{\sim}(132)$.

\item [(ii)] $Y \equiv Y^t$ for all Young diagrams $Y$.

\item [(iii)] $Y\equiv Y^{\prime}_c+Y^{\prime\prime}_c$.
\end{itemize}
\label{lem-equiv}
\end{lem}

\noindent
{\sc Proof:} Let ${\mathcal T}_{\tau}(Y)$ be the number of
transversals of $Y$ avoiding permutation $\tau$. In our previous
notation, ${\mathcal T}(Y)={\mathcal T}_{(213)}(Y)$. Since
$M(132)=M(213)^t$, ${\mathcal T}_{(132)}(Y)={\mathcal T}_{(213)}(Y^t)$.  By
definition of SWE, (i) and (ii) are equivalent.

\smallskip We will show the equivalence of (ii) and (iii) by induction
on the size $n$ of $Y$. When $n\leq 3$, (ii) and (iii) can be easily
checked by hand. Note that from the definitions of column reduction,
deletion and decomposition, $(Y^t)_c^{\prime}= (Y_r^{\prime})^t$ and
$(Y^t)_c^{\prime\prime}= (Y_r^{\prime\prime})^t$.  Assume now that
(ii) and (iii) are equivalent for size $\leq n$.

\smallskip
Assume first that (iii) is true for size $\leq n+1$. Then we can use
(ii) for sizes $\leq n$:
\[Y^t\stackrel{\text{(iii)}}{\equiv}
(Y^t)_c^{\prime}+(Y^t)_c^{\prime\prime}\stackrel{\text{def}}{=}
(Y_r^{\prime})^t+(Y_r^{\prime\prime})^t\stackrel{\text{(ii)}}{\equiv}
Y_r^{\prime}+Y_r^{\prime\prime}\stackrel{\text{RD}}{\equiv}Y.\]
This shows $Y\equiv Y^t$ and completes the proof of
$\text{(iii)}\Rightarrow \text{(ii)}$ for size $n+1$.

\smallskip 
Conversely, assume that (ii) is true for size $\leq n+1$. 
\begin{equation}
Y^t\stackrel{\text{(ii)}}{\equiv}Y\stackrel{\text{RD}}{\equiv}
Y_r^{\prime}+Y_r^{\prime\prime}\stackrel{\text{(ii)}}{\equiv}
(Y_r^{\prime})^t+(Y_r^{\prime\prime})^t\stackrel{\text{def}}{=}
(Y^t)_c^{\prime}+(Y^t)_c^{\prime\prime}.
\label{equiv}
\end{equation}
Replacing $Y$ by $Y^t$ in (\ref{equiv}), reads $Y\equiv
Y^{\prime}_c+Y^{\prime\prime}_c$.  This shows $\text{(ii)}\Rightarrow
\text{(iii)}$ for size $n+1$, and completes the proof of
$\text{(ii)}\Leftrightarrow \text{(iii)}$ and of
Lemma~\ref{lem-equiv}. \qed

\medskip
From now on, we shall refer to statement (iii) as Column-Decomposition
(CD).

\section{Commutativity of Row and Column Decompositions}
\label{section-commute}

\begin{lem} $Y\equiv Y^{\prime}_c + Y^{\prime\prime}_c$ for all
Young diagrams $Y$.
\label{lem-commute}
\end{lem}

\noindent
{\sc Proof:} Assume that the statement is true for all diagrams of
size smaller than the size of $Y$. The idea is to apply RD and the
assumed CD one after the other in different orders: this results in
representing both sides of the equality as sums of the same four
terms. Let us start with $Y$. Recall that $C_b$ and $C_t$ are the
bottom and right corners of $Y$, respectively.  Apply first \text{RD}:
\[Y\stackrel{\text{RD}}{\equiv}Y^{\prime}_r+Y^{\prime\prime}_r.\]
Next, apply the assumed CD to $Y^{\prime}_r$, and to the factor of
$Y^{\prime\prime}_r$ that still contains $C_t$, leaving the other
factor of $Y^{\prime\prime}_r$ unchanged:
\begin{equation}
Y^{\prime}_r\stackrel{\text{CD}}{\equiv}(Y^{\prime}_r)^{\prime}_c+
(Y^{\prime}_r)^{\prime\prime}_c,\,\,
Y^{\prime\prime}_r\stackrel{\text{CD}}{\equiv}
(Y^{\prime\prime}_r)^{\prime}_c+
(Y^{\prime\prime}_r)^{\prime\prime}_c.
\label{eqn-comm1}
\end{equation}
Now start with $Y^{\prime}_c+Y^{\prime\prime}_c$, and apply RD to
$Y^{\prime}_c$, and to the factor of $Y^{\prime\prime}_c$ that still
contains $C_b$, leaving the other factor of $Y^{\prime\prime}_c$
unchanged:
\begin{equation}
Y^{\prime}_c+Y^{\prime\prime}_c
\stackrel{\text{RD}}{\equiv}(Y^{\prime}_c)^{\prime}_r+
(Y^{\prime}_c)^{\prime\prime}_r+(Y^{\prime\prime}_c)^{\prime}_r+
(Y^{\prime\prime}_c)^{\prime\prime}_r.
\label{eqn-comm2}
\end{equation}
In both equations (\ref{eqn-comm1}) and (\ref{eqn-comm2}), by abuse of
notation, we wrote $(Y^{\prime\prime}_r)^{\prime}_c$,
$(Y^{\prime\prime}_r)^{\prime\prime}_c$,
$(Y^{\prime\prime}_c)^{\prime}_r$ and
$(Y^{\prime\prime}_c)^{\prime\prime}_r$ for the cross-products of
diagrams resulting from applying CD, resp. RD, to the factors
containing $C_t$, resp. $C_b$.  From the definition of row and column
deletion, it is clear that the first terms in (\ref{eqn-comm1}) and
(\ref{eqn-comm2}) are equal:
$(Y^{\prime}_r)^{\prime}_c=Y-C_b-C_t=(Y^{\prime}_c)^{\prime}_r$. We
claim that the remaining three terms also pair up as:
\begin{equation}
(Y^{\prime}_r)^{\prime\prime}_c = (Y^{\prime\prime}_c)^{\prime}_r,\,\,
(Y^{\prime\prime}_r)^{\prime}_c = (Y^{\prime}_c)^{\prime\prime}_r,\,\,
(Y^{\prime\prime}_r)^{\prime\prime}_c=
(Y^{\prime\prime}_c)^{\prime\prime}_r,
\label{eqn-pairs}
\end{equation}
except in Case IV below where
\begin{equation}
(Y^{\prime}_r)^{\prime\prime}_c \equiv (Y^{\prime}_c)^{\prime\prime}_r,\,\,
(Y^{\prime\prime}_r)^{\prime}_c \equiv (Y^{\prime\prime}_c)^{\prime}_r,\,\,
(Y^{\prime\prime}_r)^{\prime\prime}_c\equiv
(Y^{\prime\prime}_c)^{\prime\prime}_r.
\label{eqn-pairs'}
\end{equation}
Before we embark on the proofs of (\ref{eqn-pairs}-\ref{eqn-pairs'}),
note how they fit in the general outline of the proof of CD:
\begin{eqnarray*}
Y&\stackrel{\text{RD}}{\equiv}&Y^{\prime}_r+Y^{\prime\prime}_r
\,\,\text{(known)}\\ 
&\stackrel{\text{CD}}{\equiv}&
(Y^{\prime}_r)^{\prime}_c+ (Y^{\prime}_r)^{\prime\prime}_c+
(Y^{\prime\prime}_r)^{\prime}_c+
(Y^{\prime\prime}_r)^{\prime\prime}_c\,\,\text{(assumed)}\\
&\equiv&(Y^{\prime}_c)^{\prime}_r+
(Y^{\prime}_c)^{\prime\prime}_r+(Y^{\prime\prime}_c)^{\prime}_r+
(Y^{\prime\prime}_c)^{\prime\prime}_r\,\,\text{(by examining cases below)}\\
&\equiv& Y^{\prime}_c+Y^{\prime\prime}_c\,\,\text{(converse RD, known)}. \qed
\label{outline}
\end{eqnarray*}
One of the reasons that this works is that the CD-factor in
$Y^{\prime}_r$ can be anticipated from $\mathcal B$, the CD-factor in
the original $Y$.

\smallskip
The proof of (\ref{eqn-pairs}-\ref{eqn-pairs'}) depends
solely on how the row and column decomposition interact with each
other in any given Young diagram $Y$, more precisely, on the relative
position of the two $45^{\circ}$ segments used in the
decompositions. Let the {\it RD-segment} be the segment used in RD,
and let the {\it RD-factor} be the subdiagram $A_C$ determined by the
RD-segment, and similarly for CD.  Set ${\mathcal A}:=$RD-factor, and
${\mathcal B}:=$CD-factor in $Y$. To see that equations
(\ref{eqn-pairs}-\ref{eqn-pairs'}) are true, divide all Young diagrams
$Y$ into four cases: since the RD- and CD-segments are parallel to
each other, there are only four possible relative positions for
them. We will use $\stackrel{\text{sym}}{\Rightarrow}$,
$\stackrel{\text{sym}}{=}$ and $\stackrel{\text{sym}}{\equiv}$ as
shortcuts for ``by symmetry arguments''.

\smallskip
\noindent
{\bf Case I.} The RD- and CD-factors do not overlap ($\mathcal A \cap
\mathcal B=\emptyset$), i.e. the RD- and CD-segments hit $Y$'s border
before they ``come close'' to each other. Then

\begin{eqnarray*}
(Y^{\prime}_r)^{\prime\prime}_c&=&
(Y-C_b)^{\prime\prime}_c=(Y-C_b)\!\big/_{\displaystyle{\!\{(1,n),{\mathcal
B}\}}}\times {\mathcal B}=\\
&=&\left(Y\!\!\big/_{\displaystyle{\!\{(1,n),{\mathcal
B}\}}}\right)^{\prime}_r\times {\mathcal
B}=(Y^{\prime\prime}_c)^{\prime}_r;\\
\stackrel{\text{sym}}{\Rightarrow}\,\,(Y^{\prime}_c)^{\prime\prime}_r
&=&
(Y-C_t)\!\big/_{\displaystyle{\!\{(n,1),{\mathcal A}\}}}\times {\mathcal
A}{=}
(Y^{\prime\prime}_r)^{\prime}_c;\\
(Y^{\prime\prime}_r)^{\prime\prime}_c&=&
\left(Y\!\!\big/_{\displaystyle{\!\{(n,1),{\mathcal A}\}}}\times {\mathcal
A} \right)^{\prime\prime}_c= 
\left(Y\!\!\big/_{\displaystyle{\!\{(n,1),{\mathcal A}\}}}
\right)^{\prime\prime}_c \times {\mathcal A}\\
&=&
Y\!\!\big/_{\displaystyle{\!\{(n,1),(1,n),{\mathcal B},{\mathcal A}\}}}
\times {\mathcal B}\times {\mathcal A}
\stackrel{\text{sym}}{=}
(Y^{\prime\prime}_c)^{\prime\prime}_r. \qed
\end{eqnarray*}

\begin{figure}[h]
$$\psdraw{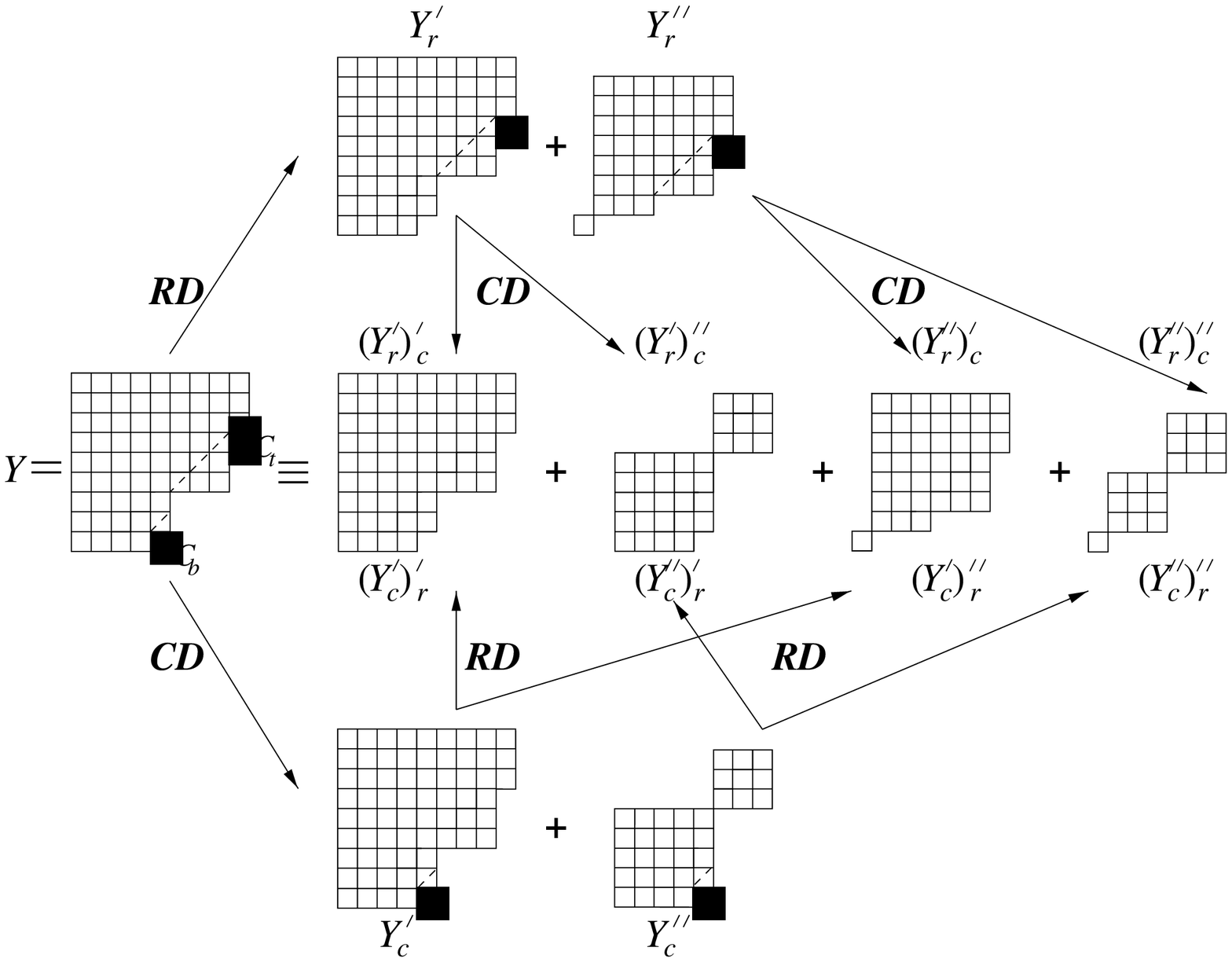}{4.9in}{3.7in}$$
\caption{Case I for $Y(9,9,9,9,8,8,5,5,5)$}
\label{comuute1-fig}
\end{figure}

\noindent
{\bf Case II.} The RD-factor contains the CD-factor ($\mathcal
A\supset \mathcal B$), i.e. the RD-segment runs ``on the inside'' of
the CD-segment.  As in Case II,
$(Y^{\prime}_r)^{\prime\prime}_c=(Y^{\prime\prime}_c)^{\prime}_r$.
Note that $C_t\in\mathcal{A}$, and therefore
$Y\!\!\big/_{\displaystyle{\!{\mathcal A}}}$ is a square.  This
justifies step ($*$) below, where ${\mathcal A}^{tr}$ denotes the top
right cell of $\mathcal A$. Note that ${\mathcal A}^{tr}$ has the same
function in $\mathcal A$ as the cell $(1,n)$ has in $Y$. The proof
works even in the extreme case where ${\mathcal A}^{tr}=(1,n)$.

\begin{eqnarray*}
(Y^{\prime}_c)^{\prime\prime}_r&=&(Y-C_t)^{\prime\prime}_r=
({\mathcal A}-C_t)\times Y\!\!\big/_{\displaystyle{\!\{(n,1),{\mathcal
A}\}}}\\
&=& {\mathcal A}^{\prime}_c \times Y\!\!\big/_{\displaystyle{\!\{(n,1),{\mathcal
A}}\}}=\left(\mathcal{A}\times Y\!\!\big/_{\displaystyle{\!\{(n,1),{\mathcal
A}}\}}\right)^{\prime}_c=(Y^{\prime\prime}_r)^{\prime}_c;\\
(Y^{\prime\prime}_r)^{\prime\prime}_c&=&
\left({\mathcal A}\times Y\!\!\big/_{\displaystyle{\!\{(n,1),{\mathcal
A}\}}}\right) ^{\prime\prime}_c={\mathcal A}^{\prime\prime}_c\times
Y\!\!\big/_{\displaystyle{\!\{(n,1),{\mathcal A}\}}}\\ &=& {\mathcal
A}\!\big/_{\displaystyle{\!\{{\mathcal B},{\mathcal A}^{tr}\}}}\times {\mathcal
B}\times Y\!\!\big/_{\displaystyle{\!\{(n,1),{\mathcal A}\}}};\\ 
&\stackrel{(*)}{=}&
\left(Y\!\!\big/_{\displaystyle{\!\{(1,n),{\mathcal
B}\}}}\right)^{\prime\prime}_r\times {\mathcal B}=
\left(Y\!\!\big/_{\displaystyle{\!\{(1,n),{\mathcal B}\}}}\times {\mathcal
B}\right)^{\prime\prime}_r=(Y^{\prime\prime}_c)^{\prime\prime}_r. \qed
\end{eqnarray*}
\begin{figure}[h]
$$\psdraw{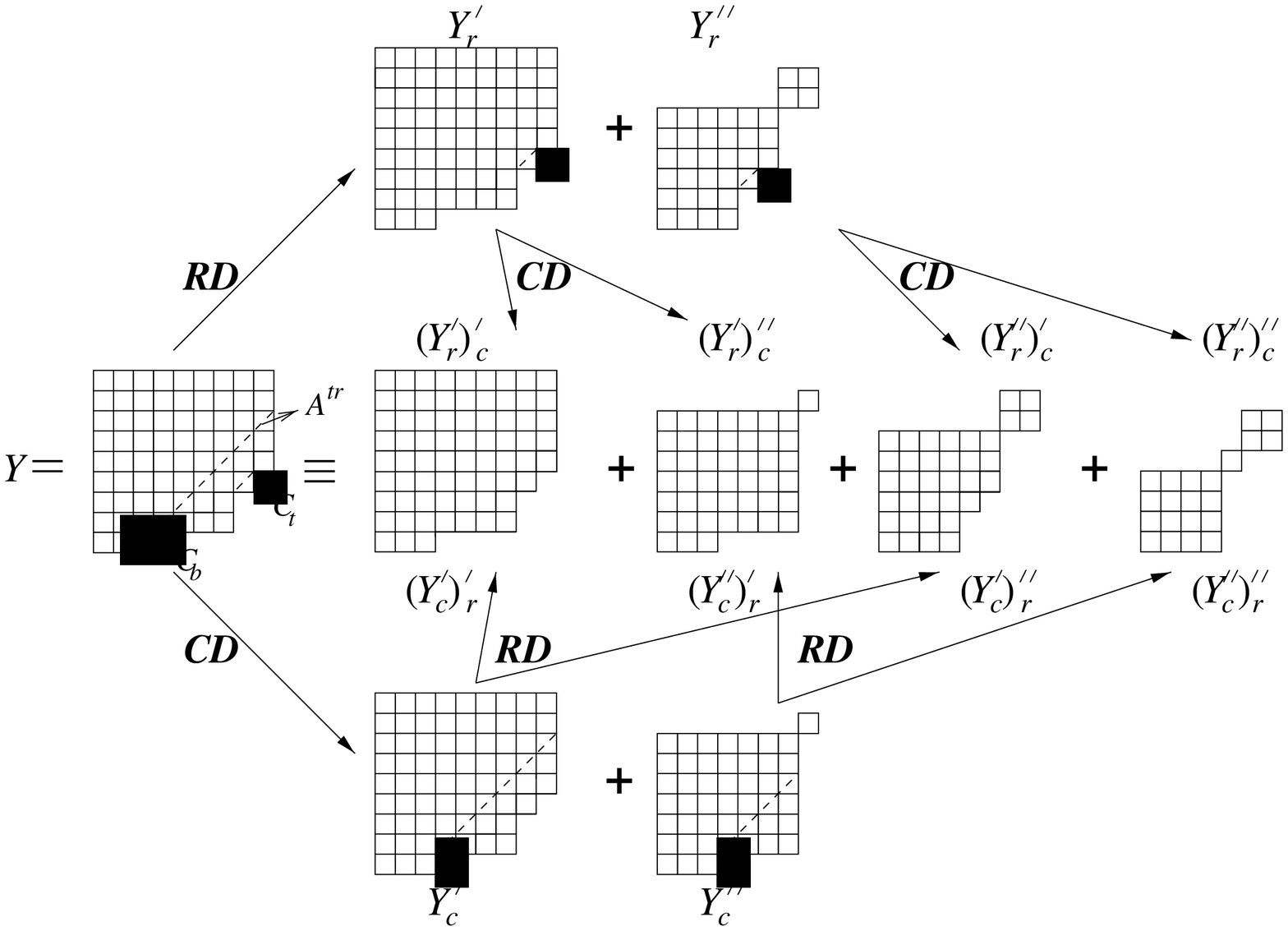}{4.9in}{3.8in}$$
\caption{Case II for $Y(9,9,9,9,9,9,7,7,4)$}
\label{commute2-fig}
\end{figure}

\noindent
{\bf Case III.} The CD-factor contains the RD-factor ($\mathcal
B\supset \mathcal A$), i.e. the
CD-segment runs ``on the inside'' of the RD-segment. This case is
symmetric to Case II. \qed

\smallskip
\noindent
{\bf Case IV.} The RD- and CD-segments overlap. This happens exactly
when the RD- and CD-segments differ from each other only in their
final cells: the RD-segment intersects the rightmost column of $Y$,
while the CD-segment intersects the bottom row of $Y$. Let ${\mathcal
D}:={\mathcal B}\!\big/_{\displaystyle{\!C_{b-1}}}= {\mathcal
A}\!\big/_{\displaystyle{\!C_{t-1}}}$.  Since ${\mathcal B}-C_b$ and
${\mathcal A}-C_t$ contain exactly one bottom, resp. rightmost, cell,
then ${\mathcal B}-C_b \equiv \mathcal{D} \equiv {\mathcal A}-C_t$ and
${\mathcal B}^{\prime\prime}_r={\mathcal D}= {\mathcal
A}^{\prime\prime}_c$.  Clearly, ${\mathcal
S}:=Y\!\!\big/_{\displaystyle{\!\mathcal D}}$ is a {\it
square}. Moreover, ${\mathcal S}_b={\mathcal
B}\!\big/_{\displaystyle{\!\mathcal D}}=C_{b-1}$, ${\mathcal
S}_t={\mathcal A}\!\big/_{\displaystyle{\!\mathcal D}}=C_{t-1}$; the
top right cell of $\mathcal S$ is the cell $(1,n)$ of $Y$, and the
bottom left cell of $\mathcal S$ is the cell $(n,1)$ of $Y$. From
this, we see that 
$Y\!\!\big/_{\displaystyle{\!\{{\mathcal
C},(1,n)\}}}= Y\!\!\big/_{\displaystyle{\!\{{\mathcal D},(n,1)\}}}
\quad \text{and} \quad Y\!\!\big/_{\displaystyle{\!\{{\mathcal
B},(1,n)\}}}= Y\!\!\big/_{\displaystyle{\!\{{\mathcal
A},(n,1)\}}}.$
Therefore,
\begin{eqnarray*}
(Y^{\prime}_r)^{\prime\prime}_c&=& (Y-C_b)^{\prime\prime}_c= {\mathcal
B}\!\big/_{\displaystyle{\!C_{b-1}}}\times
Y\!\!\big/_{\displaystyle{\!\{{\mathcal B}\big/_{C_{b-1}},(1,n)\}}}\\
&=& {\mathcal D}\times Y\!\!\big/_{\displaystyle{\!\{{\mathcal D},(1,n)\}}}
= {\mathcal D}\times
Y\!\!\big/_{\displaystyle{\!\{{\mathcal D},(n,1)\}}}
\stackrel{\text{sym}}{=} (Y^{\prime}_c)^{\prime\prime}_r;\\
(Y^{\prime\prime}_c)^{\prime}_r&=&
\left(Y\!\!\big/_{\displaystyle{\!\{{\mathcal B},(1,n)\}}}\times {\mathcal
B}\right)^{\prime}_r=
Y\!\!\big/_{\displaystyle{\!\{{\mathcal B},(1,n)\}}}\times ({\mathcal B}-C_b)\\
&\equiv& Y\!\!\big/_{\displaystyle{\!\{{\mathcal B},(1,n)\}}}\times
{\mathcal D}= Y\!\!\big/_{\displaystyle{\!\{{\mathcal
A},(n,1)\}}}\times {\mathcal D}\stackrel{\text{sym}}{\equiv}
(Y^{\prime\prime}_r)^{\prime}_c;\\
(Y^{\prime\prime}_r)^{\prime\prime}_c& = &\left({\mathcal A}\times
Y\!\!\big/_{\displaystyle{\!\{{\mathcal
A},(n,1)\}}}\right)^{\prime\prime}_c= {\mathcal
A}^{\prime\prime}_c\times Y\!\!\big/_{\displaystyle{\!\{{\mathcal
A},(n,1)\}}}=\\
&=&{\mathcal D}\times Y\!\!\big/_{\displaystyle{\!\{{\mathcal
A},(n,1)\}}}={\mathcal D}\times Y\!\!\big/_{\displaystyle{\!\{{\mathcal
B},(1,n)\}}}
\stackrel{\text{sym}}{=} (Y^{\prime\prime}_c)^{\prime\prime}_r.
\end{eqnarray*}
The three special subcases when $C_{b-1}=(n,1)$ (and hence
$C_{t-1}=(1,n)$), when $C_{b}=(n,1)$ (and hence $C_t=(1,n)$) and when
$C_b=C_t$ (i.e. $Y$ is a square), are easily checked to satisfy the
desired equalities. \qed

\smallskip
The discussion of these four cases completes the proof of Lemma 3. \qed

\vspace*{-5mm}
\begin{figure}[h]
$$\psdraw{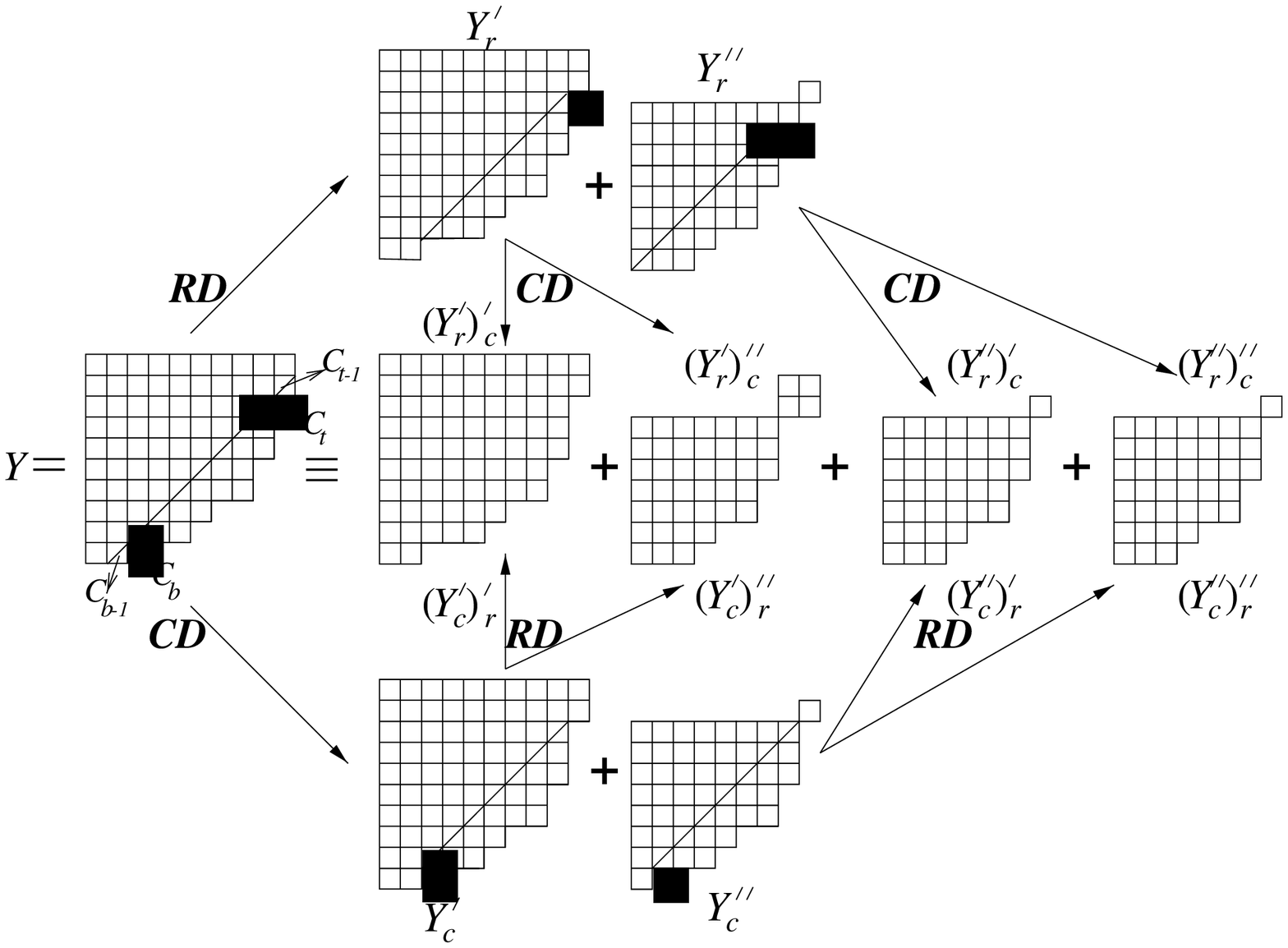}{5in}{3.7in}$$
\caption{Case IV for $Y(10,10,10,9,9,8,8,6,5,3)$}
\label{ommute3-fig}
\end{figure}
We remark that, due to the degenerate nature of Case IV, two of the
final cross-products turn out to be equal:
$(Y^{\prime\prime}_c)^{\prime}_r=
(Y^{\prime\prime}_c)^{\prime\prime}_r$, and hence
$(Y^{\prime\prime}_c)^{\prime\prime}_r$ has only two factors, rather
than the three it has in Cases I--III.

\section{New Wilf Equivalences and Consequences}
\label{section-results}

Lemmas~\ref{lem-equiv} and \ref{lem-commute} imply the main result of
the paper: $(213)\stackrel{S}{\sim}(132)$. Combined with
Proposition~\ref{prop-BW}, this establishes a new class of
Wilf-equivalent permutations.

\medskip
\noindent{\bf Theorem 1.} {\it The permutations $(213)$ and $(132)$
are shape-Wilf-equivalent. Consequently, for any $\tau\in S_{n-3}$,
the permutations $(n-1,n-2,n,\tau)$ and $(n-2,n,n-1,\tau)$ are
Wilf-equivalent.}

\smallskip
In particular, this completes the classification up to
Wilf-equivalences of $S_n$, for $n\leq 7$. Fig.~\ref{fig-summary}
lists the number of symmetry classes and Wilf-classes in each such
$S_n$.

\begin{figure}[h]
{\small \begin{tabular}{|c||c|c|c|c|c|c|c|}
\hline\hline
classes&$S_1$&$S_2$&$S_3$&$S_4$&$S_5$&$S_6$&$S_7$\\\hline
\hline
symmetry& 1 & 1 & 2 & 7 & 23 & 115 & 694\\ \hline
Wilf    & 1 & 1 & 1 & 3 & 16 &  91 & 595\\
\hline
\end{tabular}}
\caption{Number of symmetry vs. Wilf classes in $S_n$, $n\leq 7$}
\label{fig-summary}
\end{figure}

\smallskip
An amusing corollary about numerical equivalence of Young diagrams can
be deduced from the above theorem and the row-decomposition formula.
Recall that $St$ is the standard staircase diagram.  The {\it
$k$-staircase} $St_k$ is the Young diagram which consists of $St$ plus
the full $k-1$ diagonals below the diagonal of $St$. In particular,
$St$ is $St_1$, and the square $n\times n$ is $St_n$. The {\it
critical} staircase $St_k$ of $Y$ is the first staircase whose
complement $Y\backslash St_k$ is a union of at least two connected
components. Label such components by $St_k^j$, for $j=1,2,...$,
starting at the bottom left corner of $Y$.  Thus, for every $Y$ with
critical staircase $St_k$, we have the {\it critical} decomposition
$Y=St_k \cup_j St_k^j$.

\begin{cor} Let $Y=St_k \cup_{j=1}^l St_k^j$ be the critical
decomposition of $Y$. The operations of permuting and of transposing
the components $St_k^j$ result in Young diagrams numerically
equivalent to $Y$.
\label{cor-hanging}
\end{cor}

\begin{figure}[h]
$$\psdraw{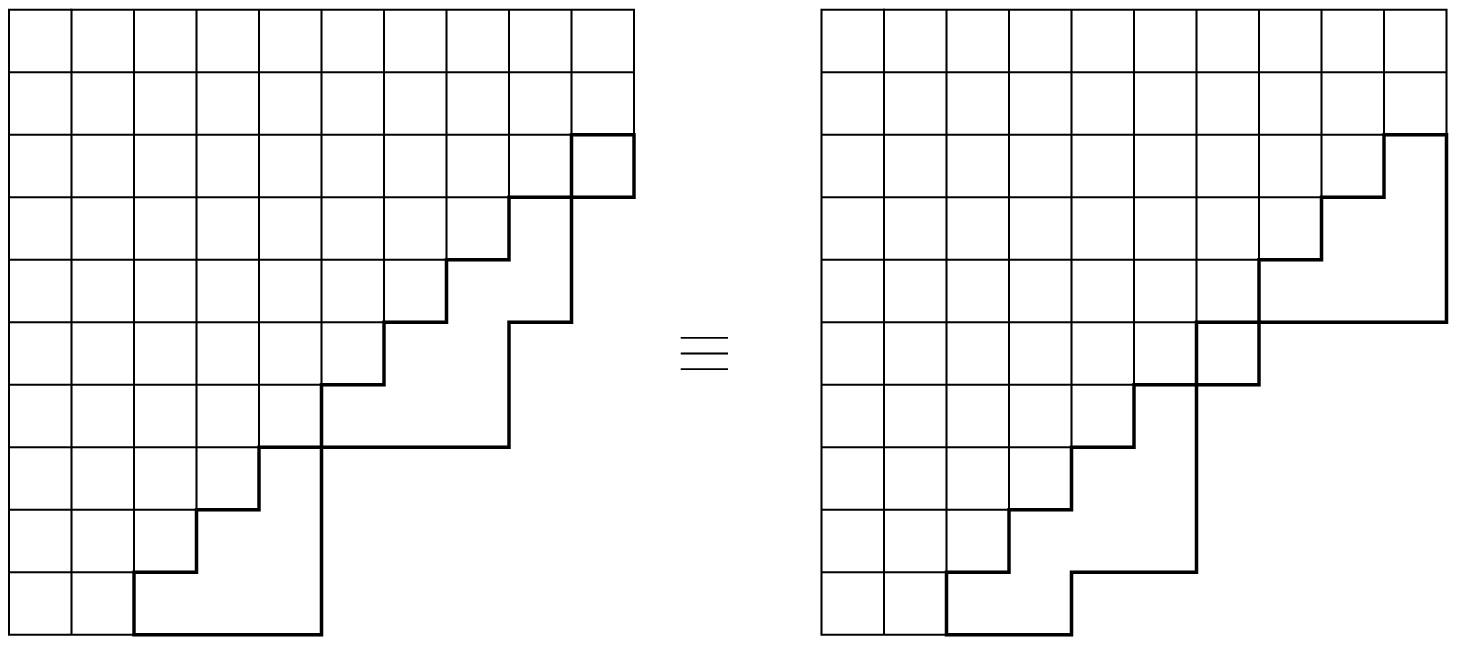}{2in}{0.9in}$$
\caption{$St_2\cup St_2^1\cup St_2^2\cup St_2^3
\equiv St_2\cup (St_2^2)^t\cup St_2^3\cup St_2^1$}
\label{hanging-fig}
\end{figure}

In other words, let $\tau\in S_l$ be any permutation, and let
$\vec{t}=(t_1,t_2,...,t_l)\in \{1,t\}^l$ correspond to a choice $t_j=t$ to
transpose, resp. $t_j=1$ not to transpose, the component
$St_k^{\tau(j)}$. Then the following (ordered) critical decompositions
represent numerically equivalent Young diagrams
(cf. Fig.~\ref{hanging-fig}):
\[St_k\bigcup_{j=1}^l St_k^j \equiv St_k \bigcup_{j=1}^l 
\big(St_k^{\tau(j)}\big)^{t_j}.\]

\noindent
{\sc Proof:} We use induction on the size of $Y$. The initial cases
are easily verified. Moreover, when there is only {\it one}
``hanging'' shape, the corollary simply states that transposing $Y$
will yield a numerically equivalent diagram $Y^t$: this is the content
of Theorem~\ref{thm-Wilf}.

Suppose now that there are at least two hanging shapes: $St_k^1$ is
the bottom shape, and let's name the remaining shapes the ``upper''
shapes.  Choose a permutation $\sigma$ and a transposition vector
$\vec{t}$ both of which leave the bottom shape $St_k^1$ fixed, and
operate on the upper shapes of $Y$.  Apply the row-decomposition
formula to $Y$: when $C$ runs over all bottom cells of $Y$, we have
\begin{equation}
Y\equiv \sum_C Y_C=\sum_C A_C\times B_C.
\label{pf-cor}
\end{equation}
From the definition of the critical decomposition of $Y$, all
RD-factors $A_C$ are parts of the bottom shape $St_k^1$, except for
the first and the last RD-factors: there $A_C=Y_r$ or $A_C=\emptyset$.  In
any case, none of the upper shapes $St_k^2,\,\,St_k^3,...$ are
broken up in (\ref{pf-cor}). Thus, by induction hypothesis, we can
apply $\sigma$ and $\vec{t}$ to each $Y_C$:
$\big(Y_C^{\sigma}\big)^{\vec{t}}= A_C\times
\big(B_C^{\sigma}\big)^{\vec{t}}$. We can then put back together the
resulting diagrams via another application of the row-decomposition
formula:
\begin{equation*}
Y\equiv \sum_C Y_C=\sum_C A_C\times B_C \stackrel{\text{ind.}}{\equiv}
\sum_C A_C\times \big(B_C^{\sigma}\big)^{\vec{t}} \equiv
\big(Y^{\sigma}\big)^{\vec{t}}.
\end{equation*}
This shows that leaving the bottom shape fixed, we can permute and
transpose the upper shapes in any way we like. We also know, again
from Theorem~\ref{thm-Wilf}, that transposing the whole diagram $Y$
yields a numerically equivalent diagram $Y^t$. It is an easy exercise
in algebra to verify that these two types of operations generate the
whole group of operations required in
Corollary~\ref{cor-hanging}. \qed

\smallskip
The conclusion of Corollary~\ref{cor-hanging} holds under a slightly
relaxed hypothesis regarding which staircases can be used instead of the
critical staircase: as long as the RD (or CD) formula breaks up only
the bottom (or only the top) hanging shape, the above proof goes
through without modifications. Further generalizations are also
possible, for instance, applying recursively the Corollary just within
a hanging shape. Finally, this can all be used to write down a
generating function for the numbers $T_{(213)}(Y)$, but we will not do
this here since it will take us too far afield.

\section{Further Discussion}
\label{section-conjectures}

A careful investigation of the new Wilf-pair $(546213)\sim(465213)$ in
$S_6$, leads to the observation that both permutation matrices can be
decomposed into two blocks of $3\times 3$ matrices, and further, that
moving from one decomposition to the other involves a transposition of
one of the blocks. Thus, one might be lead to conjecture
that for any permutation matrices $A$ and $B$, the
following permutation matrices are Wilf-equivalent:
\begin{equation}
\left( \begin{array}{c|c}
          A & 0\\ \hline
          0 & B
      \end{array} \right) {\sim}
\left(\begin{array}{c|c}
          A^t & 0\\ \hline
          0 & B
      \end{array} \right)\cdot
\label{conj-strong}
\end{equation}

\noindent In order for (\ref{conj-strong}) to give any new
Wilf-classes, other than those obtained by symmetry or $I_t
\stackrel{s}{\sim}J_t$, both $A$ and $B$ must be non-symmetric
matrices. In $S_6$, there is only one such pair up to symmetry (denote
the 1's by dots, and omit all 0's):
\[{\tiny M(546213)=\left( \begin{array}{ccc|ccc}
          & & \!\!\bullet\!\! & &&\\
          \!\!\bullet\!\!&& & &&\\
          & \!\!\bullet\!\! &&&\\ \hline
          &&&& & \!\!\bullet\!\! \\
          &&&\!\!\bullet\!\!&& \\
          &&& &\!\!\bullet\!\! &\\ 
      \end{array} \right) {\sim}
     \left(\begin{array}{ccc|ccc}
          & \!\!\bullet\!\! &&&\\  
          & & \!\!\bullet\!\! & &&\\
           \!\!\bullet\!\!&& & &&\\\hline
          &&&& & \!\!\bullet\!\! \\
          &&&\!\!\bullet\!\!&& \\
          &&& &\!\!\bullet\!\! &\\ 
      \end{array} \right)=M(465213)\cdot}\] 
\noindent 
In $S_7$, there are essentially 7 new Wilf-pairs which are covered by
(\ref{conj-strong}). Not surprisingly, these are the same
pairs appearing in Fig.~\ref{new-Wilf-S7}. One possible approach to
prove (\ref{conj-strong}) would be to show
$A\stackrel{s}{\sim} A^t$ for any permutation matrix
$A$. Unfortunately, this is {\it not} true; it fails already in $S_4$,
e.g.  $(3142)\not\stackrel{s}{\sim}(2413)$ since
\[|S_{(6,6,6,6,5,5)}(3142)|=394< 395=|S_{(6,6,6,6,5,5)}(2413)|.\] The
Wilf-equivalence in (\ref{conj-strong}) also {\it fails} in $S_8$. For
example:
\[{\tiny  M(68572413)=
       \left( \begin{array}{cccc|cccc}
          & \!\!\bullet\!\!& &&& &&\\
          &&&\!\!\bullet\!\!&&&&\\
          \!\!\bullet\!\!&&&&&&&\\
          & &\!\!\bullet\!\! &&&&\\ \hline
          &&&&& \!\!\bullet\!\!&&\\
          &&&&&&&\!\!\bullet\!\!\\
          &&&&\!\!\bullet\!\!&&&\\
          &&&&& &\!\!\bullet\!\!& \\
      \end{array} \right) {\not\sim}
     \left(\begin{array}{cccc|cccc}
          & &\!\!\bullet\!\! &&&&\\ 
          \!\!\bullet\!\! &&&&&&\\ 
          &&&\!\!\bullet\!\!&&&&\\
          & \!\!\bullet\!\!& &&& &&\\\hline
          &&&&& \!\!\bullet\!\!&&\\
          &&&&&&&\!\!\bullet\!\!\\
          &&&&\!\!\bullet\!\!&&&\\
          &&&&& &\!\!\bullet\!\!& \\
      \end{array} \right)=M(75862413)\cdot}\] 
The two permutations are equinumerant up to level 11, but
they split on level 12:
\[|S_{12}(75862413)|=476576750<476576751=
|S_{12}(68572413)|.\] 
This forces a reexamination of the 8 new
Wilf-pairs in $S_6$ and $S_7$. If we choose the ``right''
representatives of the symmetry classes, we can see that each
permutation matrix contains a block corresponding to $(213)$ or
$(132)$. This led to conjecturing and proving the
shape-Wilf-equivalence $(213)\stackrel{s}{\sim}(132)$.  As we can see,
this SWE is far from coincidental, and it is the reason for the
infinitely many new Wilf-equivalences of Theorem~\ref{thm-Wilf}.

\medskip
Let us now shift the emphasis of our discussion to a slightly
different question. If we write down a table enumerating $|S_n(\tau)|$
for a fixed $\tau\in S_4$ as $n$ increases, we notice a very plausible
conjecture (cf. \cite{We0}):

\begin{figure}[h]
{\small \begin{tabular}{|c||c|c|c|c|c|c|c|c|}
\hline\hline
$\tau\in S_4$&$S_6(\tau)$&$S_7(\tau)$
&$S_8(\tau)$&$S_9(\tau)$&$S_{10}(\tau)$&$S_{11}(\tau)$&$S_{12}(\tau)$
&$S_{13}(\tau)$\\\hline
\hline
(4132) & 512 & 2740 & 15485 & 91245 & 555662 & 3475090 &
22214707 & 144640291\\ \hline
(1234) & 513 & 2761 & 15767 & 94359 & 586590 & 3763290 &
24792705 & 167078577\\\hline
(1324) & 513 & 2762 & 15793 & 94776 & 591950 & 3824112 & 
25431452 & 173453058\\\hline
\end{tabular}}
\caption{Classification of $S_4$ up to Wilf-equivalence}
\label{fig-S4}
\end{figure}

\begin{conj}
If $|S_k(\tau)| < |S_k(\sigma)|$ for some $k$, then $|S_n(\tau)| <
|S_n(\sigma)|$ for all $n\geq k$. In other words, modulo Wilf-equivalence,
we can order all permutations in $S_n$ according to their relative
restrictiveness: $\tau<\sigma$ if $|S_k(\tau)|<|S_k(\sigma)|$ for some
$k$.
\label{strong-order}
\end{conj}

$S_4$ is the first non-trivial case of Conjecture~\ref{strong-order}.
It was partially proved by Bona in \cite{Bon1}, where he shows
$|S_n(1423)|<|S_n(1324)|$ for $n\geq 6$ and $|S_n(1234)|<|S_n(1324)|$
for $n\geq 7$, in relation to a conjecture of Wilf and Stanley. To the
best of our knowledge, no one had published a counterexample to
Conjecture~\ref{strong-order}, until we found a counterexample in
$S_5$, followed by various types of counterexamples in $S_6$ and
$S_7$.

\begin{figure}[h]
{\tiny \begin{tabular}{|c||c|c|c|c|c|c|l|}
\hline\hline
$\tau\in S_5$&$S_7(\tau)$&$S_8(\tau)$&$S_9(\tau)$
& $S_{10}(\tau)$ & $S_{11}(\tau)$ & $S_{12}(\tau)$ & $S_{13}(\tau)$\\\hline
\hline
(25314)& 4578 & 33184 & 258757 & 2136978 & 18478134 & 165857600 & 1535336290\\ 
\hline
(31524)& 4579 & 33216 & 259401 & 2147525 & 18632512 & 167969934 & 1563027614\\ 
\hline
(35214)& 4579 & 33218 & 259483 & 2149558 & 18672277 & 168648090 & 1573625606\\ 
\hline
(35124)& 4580 & 33249 & 260092 & 2159381 & 18815124 & 170605392 & 1599499163\\ 
\hline
(53124)& 4580 & 33252 & 260202 & 2161837 & 18858720 & 171285237 & 1609282391\\ 
\hline
(42351)& 4580 & 33252 & 260204 & 2161930 & 18861307 & 171341565 & 1610345257\\ 
\hline
(35241)& 4580 & 33254 & 260285 & 2163930 & 18900534 & 172016256 & 1621031261\\ 
\hline
(53241)& 4580 & 33256 & 260370 & 2166120 & 18945144 & 172810050 & 1633997788*\\\hline
(43251)& 4581 & 33283 & 260805 & 2171393 & 18994464 & 173094540 & 1632480259*\\ \hline
(32541)& 4581 & 33284 & 260847 & 2172454 & 19015582 & 173461305 & 1638327423\\ 
\hline
(34215)& 4581 & 33285 & 260886 & 2173374 & 19032746 & 173741467 & 1642533692\\ 
\hline
(31245)& 4581 & 33286 & 260927 & 2174398 & 19053058 & 174094868 & 1648198050\\ 
\hline
(42315)& 4581 & 33287 & 260967 & 2175379 & 19072271 & 174426353 & 1653484169\\ 
\hline
(12345)& 4582 & 33324 & 261808 & 2190688 & 19318688 & 178108704 & 1705985883\\ 
\hline
(53421)& 4582 & 33325 & 261853 & 2191902 & 19344408 & 178582940 & 1713999264\\ 
\hline
(52341)& 4582 & 33325 & 261863 & 2192390 & 19358590 & 178904675 & 1720317763\\ 
\hline
\end{tabular}}
\caption{Classification of $S_5$ up to Wilf equivalence}
\label{fig-S5}
\end{figure}

As Fig.~\ref{fig-S5} suggests, $(53241)$ and $(43251)$ {\it cannot} be
ordered since $S_k(53241)<S_k(43251)$ for $k\leq 12$, but
$S_{13}(53241)>S_{13}(43251)$. Fig.~\ref{counter-S6}-\ref{counter-S7}
list all counterexamples up to level $13$ in $S_6$, and some
counterexamples in $S_7$. The asterisks indicate the first level at
which the corresponding permutations ``switch'' their relative
restrictiveness, and hence cannot be ordered as in
Conjecture~\ref{strong-order}. The ``!!'' in Fig.~\ref{counter-S6}
refers to the permutation (546213), which is also part of the new
Wilf-equivalence of Theorem~\ref{thm-Wilf}.

\begin{figure}[h]
{\tiny \begin{tabular}{|c||c|l|l|l|l|}
\hline\hline
$\tau\in S_7$&$S_9(\tau)$&$S_{10}(\tau)$&$S_{11}(\tau)$&$S_{12}(\tau)$ &
$S_{13}(\tau)$\\\hline
\hline
253146 & 39424 & 344580 & 3283521 & 33633237* & 366084190 \\    
523614 & 39425 & 344611 & 3283955 & 33632674* & 365858205 \\\hline
352164 & 39425 & 344619 & 3284418 & 33648781* & 366298292 \\    
426153 & 39425 & 344620 & 3284441 & 33648549* & 366268369 \\\hline
365241 & 39425 & 344633 & 3285228 & 33676816* & 367058930 \\    
236145 & 39426 & 344661 & 3285505 & 33671423* & 366717782 \\\hline
356214 & 39426 & 344671 & 3286054 & 33689584  & 367192027* \\   
315642 & 39426 & 344672 & 3286086 & 33689894  & 367181171* \\\hline        
315426 & 39426 & 344678 & 3286493 & 33705839* & 367659357 \\    
524316 & 39426 & 344679 & 3286521 & 33705830* & 367633865 \\\hline
325146 & 39426 & 344681 & 3286660 & 33711372* & 367802282 \\    
532416 & 39426 & 344682 & 3286686 & 33711239* & 367772447 \\\hline
326154 & 39427 & 344724 & 3287748 & 33732379* & 368139497 \\    
156324 & 39427 & 344724 & 3287749 & 33732314* & 368131670 \\\hline
643251 & 39427 & 344725 & 3287840 & 33736582* & 368277055 \\    
453261 & 39427 & 344726 & 3287851 & 33735408* & 368204877 \\\hline
435216 & 39427 & 344726 & 3287877 & 33737112  & 368270512* \\     
532164 & 39427 & 344727 & 3287904 & 33737204  & 368255863* \\\hline
632541 & 39427 & 344727 & 3287974*& 33741840 & 368436494 \\   
513246 & 39427 & 344728 & 3287971*& 33739711 & 368325636 \\\hline
624351 & 39427 & 344732 & 3288259 & 33751604  & 368699773* \\     
652413 & 39427 & 344733 & 3288292 & 33751853  & 368681504* \\\hline
643521 & 39427 & 344735 & 3288426 & 33757116  & 368841108* \\     
546213 & 39428 & 344772 & 3289163 & 33765743  & 368833207*!! \\\hline\hline
\end{tabular}}
\caption{Counterexamples to Conjecture~\ref{strong-order} in
$S_6$}
\label{counter-S6}
\end{figure}

\begin{figure}[h]
{\tiny \begin{tabular}{|c||c|l|l|l|l|}
\hline\hline
$\tau\in S_7$&$S_9(\tau)$&$S_{10}(\tau)$&$S_{11}(\tau)$&$S_{12}(\tau)$&
$S_{13}(\tau)$\\\hline\hline
4725631	& 361297 & 3587548 & 38941876** & 457096663* & 5746716800\\
5427316	& 361297 & 3587549 & 38941876** & 457093463* & 5746500592\\\hline
7536241	& 361297 & 3587572 & 38943767*  & 457184933  & 5749921907 \\
2645173	& 361298 & 3587605 & 38943703*  & 457138045  & 5747366258\\\hline
\end{tabular}}
\caption{Some counterexamples to Conjecture~\ref{strong-order} in
$S_7$}
\label{counter-S7}
\end{figure}

\smallskip
However, the above counterexamples do not preclude an {\it asymptotic}
ordering of all permutations in $S_n$.

\begin{defn}{\rm 
For $\tau,\sigma\in S_n$, we say that $\tau$ is {\it asymptotically
smaller} than $\sigma$ if $|S_k(\tau)| < |S_k(\sigma)|$ for all
$k\gg 0$.
}\end{defn}

\begin{conj} 
We can order asymptotically all permutations
in $S_n$, modulo Wilf-equivalence.
\label{weak-order}
\end{conj}
\noindent Regev \cite{Re} and Bona \cite{Bon2} have worked on
asymptotic behavior of certain types of permutations, but as of now,
Conjecture~\ref{weak-order} and some possible modifications of it are
far from proven.

\section* {Acknowledgments}

Zvezdelina Stankova-Frenkel is grateful to Bernd Sturmfels (UC
Berkeley) for his helpful suggestions both on mathematical and
computational issues related to the project; to David Moews
(University of Connecticut) for writing two of the computer programs
used in this project; and to Paulo de Souza (UC Berkeley) and Tom
Davis (Palo Alto) for their generous help in installing, testing and
running the necessary computer software.  Both authors would like to
thank Olivier Guibert (Universit\'e Bordeaux 1, Talence Cedex, France)
for providing his computer program on enumeration of forbidden
subsequences.

\end{document}